\documentclass[a4paper,fleqn]{cas-sc}
\usepackage{amsmath}
\usepackage{ifpdf}
\usepackage{indentfirst}
\usepackage{hyperref}
\usepackage{algorithm}
\usepackage{algorithmicx}  
\usepackage{algpseudocode}
\usepackage{subcaption}
\usepackage{booktabs}
\usepackage{multirow}
\usepackage{float}
\usepackage{enumitem}
\usepackage{adjustbox}

\numberwithin{equation}{section}

\usepackage{graphicx,cite,amssymb,lineno,color,amsfonts}

\newtheorem{Theorem}{Theorem}[section]
\newtheorem{Lemma}{Lemma}[section]

\newtheorem{Definition}{Definition}[section]
\newtheorem{Remark}{Remark}[section]

\newtheorem{Assumption}{Assumption}[section]

\makeatletter
\def\elsartstyle{%
    \def\normalsize{\@setfontsize\normalsize\@xiipt{14.5}}
    \def\small{\@setfontsize\small\@xipt{13.6}}
    \let\footnotesize=\small
    \def\large{\@setfontsize\large\@xivpt{18}}
    \def\Large{\@setfontsize\Large\@xviipt{22}}
    \skip\@mpfootins = 18\p@ \@plus 2\p@
    \normalsize
} \@ifundefined{square}{}{\let\Box\square} \makeatother

\def\file#1{\texttt{#1}}

\begin{document}
\let\WriteBookmarks\relax
\let\printorcid\relax 

\shorttitle{}    
\shortauthors{}  
\title [mode = title]{An inexact semismooth Newton-Krylov method for semilinear elliptic optimal control problem}  

\tnotemark[1] 


\tnotetext[1]{This work was supported by the Guangdong Basic and Applied Basic Research Foundation of China.}

\author[1]{Shiqi Chen}

\credit{}



\author[1]{Xuesong Chen}
\cormark[1]
\ead{chenxs@gdut.edu.cn}

\credit{}

\affiliation[1]{organization={School of Mathematics and Statistics, Guangdong University of Technology},
            city={Guangzhou},
            postcode={510520}, 
            country={P.R. China}}
\cortext[1]{Corresponding author}


\begin{abstract}
  An inexact semismooth Newton method has been proposed for solving semi-linear elliptic optimal control problems in this paper. This method incorporates the generalized minimal residual (GMRES) method, a type of Krylov subspace method, to solve the Newton equations and utilizes nonmonotonic line search to adjust the iteration step size. The original problem is reformulated into a nonlinear equation through variational inequality principles and discretized using a second-order finite difference scheme. By leveraging slanting differentiability, the algorithm constructs semismooth Newton directions and employs GMRES method to inexactly solve the Newton equations, significantly reducing computational overhead. A dynamic nonmonotonic line search strategy is introduced to adjust stepsizes adaptively, ensuring global convergence while overcoming local stagnation. Theoretical analysis demonstrates that the algorithm achieves superlinear convergence near optimal solutions when the residual control parameter $\eta_k$ approaches to 0. Numerical experiments validate the method's accuracy and efficiency in solving semilinear elliptic optimal control problems, corroborating theoretical insights.
\end{abstract}




\begin{keywords}
  Optimal control\sep 
  Inexact semismooth Newton method\sep 
  Krylov subspace method\sep
  Superlinear convergence\sep
  Global convergence\sep
\end{keywords}

\maketitle

\section{Introduction}
We focus on an optimal control problem governed by a semilinear elliptic partial differential equation (PDE), subject to box constraints, as detailed subsequently:
\begin{equation}
  \label{equation 1.1}
  \min_{u}J(u)=\frac{1}{2}{\Vert  y-y_d \Vert ^2_{L^2(\Omega)}}+\frac{\alpha}{2}{\Vert  u \Vert ^2_{L^2(\Omega)}}
\end{equation}

\begin{equation}
  \label{equation 1.2}
   s.t. \left\{
  \begin{array}{rl}
  -\varDelta y+S(y)=f+u & \quad \text{in}\;\Omega,\\
  y = 0 & \quad \text{on}\;{\partial}{\Omega},\\
  u_a{\le}u{\le}u_b & \quad \text{a.e.\;in} \;\Omega,
  \end{array}\right.
\end{equation}
where $\varDelta$ is Laplace operator, $\Omega = (0,1)^2$. $u \in L^2(\Omega)$ is control variable, $y\in H^1_0(\Omega)$ is state variable. $(y_d,f) \in (H^1_0(\Omega),L^2(\Omega))$ are given functions.
$S:R \rightarrow R $ is a nonlinear function. $ \alpha>0 $ is a positive constant. $ \{u_a,u_b\} \subset L^{\infty}(\Omega)$, $u_a\leq u \leq u_b$ is box constraint.

The numerical methods of the optimal control problem with constraints of PDE are mainly divided into two distinct strategies: the discretize-then-optimize (DO) approach and the optimize-then-discretize (OD) approach. Liu \cite{MR3987421} examines an elliptic distributed optimal control problem and compares the convergence of the OD and DO methods. There is no absolute superiority of one method over the other; the most effective approach depends on the suitability of the approach for the specific problem at hand. In this paper, the OD approach is adopted.
 
For the OD approach, it begins by discussing the first or second order necessary continuous optimal conditions of the original problem, that is, the equations composed of the equation of state, the costate variable state and the variational inequality, and then the discretization of the continuous optimal system is carried out by an appropriate discretize scheme, such as finite difference method and finite element method, the obtained discrete linear/nonlinear system can be solved by various well-designed efficient iterative algorithm. For example, in \cite{MR4419342}, the authors derive the optimality conditions for an elliptic distributed optimal control problem using the calculus of variations. They then discretize these conditions using a multiscale finite element method and finally solve the resulting discrete system to obtain numerical solutions. In another paper \cite{MR3817765} by the authors, they employ $C_0$ interior penalty methods to discretize the optimality conditions of the elliptic distributed optimal control problem with pointwise state constraints and subsequently solve for the numerical solution. There are plenty of numerical methods are employed to slove (\ref{equation 1.1})-(\ref{equation 1.2})\cite{MR4447697,MR250665,MR3461397,MR3768951,MR3374552,MR4636154}. Zhang \cite{MR4226251} proposed a conjugate gradient algorithm to solve distributed elliptic optimal control problems and prove the algorithm has global convergence property and linear convergence rate. In \cite{MR4174643}, authors employ an augmented Lagrangian method to address a class of semilinear elliptic optimal control problems with pointwise state constraints. In addition, stationary points of the subproblems are shown to exist within arbitrarily small neighborhoods of local solution with an additional quadratic growth condition. 

The Newton method is known for its rapid convergence property. Smoothing methods and Semismooth methods had been studied for semismooth problems\cite{MR1458218,MR1443628,MR1786137,MR1734657}. Chen \cite{MR1786137} introduced the concept of slant differentiability and developed a unified theoretical framework to investigate the convergence properties of smoothing methods and semismooth methods. The analysis demonstrated that these methods achieve superlinear convergence rates. To enhance the convergence rate of the semismooth Newton method, Liu utilize a multigrid method in \cite{MR3461397}. 

Numerous semismooth Newton methods with global convergence properties have been proposed. \cite{MR2278446}Global convergence can be achieved through various numerical optimization techniques, including but not limited to line search methods, dogleg methods and trust region methods. Zhang \cite{MR4447697} brought up a semismooth Newton method for semilinear optimal control problem with $L_1$ control cost. To address the limitations of local convergence of the Newton method, Zhang introduced line search techniques to ensure global convergence of the algorithm. 

When employing (semi)smooth Newton methods to solve linear system (Newton-step equation) for obtaining descent direction, the computational cost becomes prohibitive for large-scale matrices due to excessive memory requirements and prolonged computation time. This computational challenge motivates the development of inexact solving strategies, where the Newton-step equation is solved only to a prescribed residual tolerance rather than exact precision \cite{MR1119247,MR2589588,MR1273766,  MR4891821,MR4925867}. Compared with exact solves for descent direction, inexact approache circumvents the numerical challenges associated with large-scale matrix inversion.

Direct solvers face significant challenges when solving large-scale linear systems due to their high computational cost and memory requirements, particularly for sparse or high-dimensional problems. Krylov subspace methods, such as the Generalized Minimal Residual (GMRES) method, overcome these limitations by leveraging iterative approximation within low-dimensional subspaces. Rather than explicitly factorizing the matrix, these methods rely on efficient matrix-vector multiplications, drastically reducing memory usage and computational overhead. We can find the application of this method in \cite{MR2589588,MR1860971,MR2292433}.

The main contributions of this work are threefold. First, we derive the first-order necessary optimality conditions for the original problem (\ref{equation 1.1})-(\ref{equation 1.2}) and discretize them using a five-point difference scheme, leading to a nonlinear system with challenging nonsmooth terms. To address this, we develop a novel solution framework incorporating slant differentiability concepts and a semismooth Newton method. Second, we design an efficient inexact solution strategy where the Newton equations are solved approximately using GMRES, combined with a nonmonotone backtracking line search to ensure global convergence. Our analysis establishes the local superlinear convergence rate of this approach. Finally, comprehensive numerical experiments demonstrate the computational efficiency and precision of the proposed method, while validating the theoretical convergence results.

The remainder of this paper is organized as follows. Section 2 establishes the first-order necessary optimality conditions for problem (\ref{equation 1.1})-(\ref{equation 1.2}), which are then discretized using a five-point difference scheme to yield a solvable nonlinear system. In Section 3, we develop the Inexact Newton-GMRES method with nonmonotone line search for solving the resulting nonlinear system. The feasibility and convergence analysis of this method are rigorously examined in Section 4. Finally, Section 5 presents numerical experiments to validate the theoretical findings.

\section{First-order optimality system and its numerical discretization}

In this section, we derive the first-order necessary optimality conditions for problem (\ref{equation 1.1})-(\ref{equation 1.2}). Due to the non-convex nature of the problem, uniqueness of the optimal solution cannot be guaranteed a priori. To establish solution uniqueness, we first introduce the following key assumptions.
\cite{MR250665}Assume that the nonlinear term $S(y)$ satisfies the condition: $S(\cdot)\in W^{2,\infty}(-R,R)$ for any $R>0$,
$S'(y)\in L^2(\Omega)$ for any $y \in H^1(\Omega)$, and $S'\geq 0$. Let 
\begin{equation*}
  \begin{gathered}
    a(\bar{x},\bar{y}) = \int_{\Omega} \nabla \bar{x} \nabla \bar{y} \, \mathrm{d}\Omega ,\\
    (\bar{x},\bar{y}) = \int_{\Omega}\bar{x} \bar{y} \, \mathrm{d}\Omega
  \end{gathered}
\end{equation*}
for each $\bar{x},\bar{y}\in L^2(\Omega)$.
Thus the weak formulation
\begin{equation*}
  a(y,w)+(S(y),w)-(u,w)=(f,w), \quad \forall w \in H^1_0(\Omega)
\end{equation*}
has a unique solution. Meanwhile, since (\ref{equation 1.1}) is convex,
\begin{equation*}
  U_{ad}= \{u_a{\le}u{\le}u_b\}
\end{equation*}
also is convex. Therefore, 
\begin{equation}
  \label{equation 2.1}
  \begin{gathered}
    \min_{u}J(u)=\frac{1}{2}{\Vert  y-y_d \Vert ^2_{L^2(\Omega)}}+\frac{\alpha}{2}{\Vert  u \Vert ^2_{L^2(\Omega)}} \\
    s.t. \left\{
    \begin{array}{cl}
      a(y,w)+(S(y),w)-(u,w)=(f,w) & \quad \forall w \in H^1_0(\Omega),\\
      u_a{\le}u{\le}u_b & \quad a.e.in \;\Omega
    \end{array}\right.
  \end{gathered}
\end{equation}
possesse a unique optimal solution. To describe the optimal solution of problem (\ref{equation 2.1}), the first-order necessary optimal conditions
are stated as follows in \cite{troltzsch2010optimal}:
\begin{equation}
  \label{equation 2.2}
  \begin{gathered}
    -\varDelta y+S(y)-u=f  \quad \text{in}\;\Omega \quad \text{and} \quad y = 0  \quad \text{on}\;{\partial}{\Omega}, \\
    -\varDelta p+S'(y)p+y=y_d  \quad \text{in}\;\Omega \quad \text{and} \quad p = 0  \quad \text{on}\;{\partial}{\Omega}, \\
    (\alpha u-p, v-u) \geq 0 \quad \text{for} \; \text{all} \; v \in U_{ad}.
  \end{gathered}
\end{equation}
We present the Lemma \ref{Lemma2.1} to address problem (\ref{equation 2.2}).

\begin{Lemma}
  \label{Lemma2.1}
    Using the variational inequality in (\ref{equation 2.2}), we get an equivalent way to describe the optimal control:
    \begin{equation}\label{equation 2.3}
      u^*=\Phi\left(\frac{p}{\alpha}\right)=\max{\left\{u_a,\min\left\{u_b,\frac{p}{\alpha}\right\}\right\}}.
    \end{equation}
\end{Lemma}

\noindent{\it Proof.} ~~ \cite{troltzsch2010optimal}The necessary and sufficient condition for variational inequality in (\ref{equation 2.2}) is given by
\begin{equation}
  \label{equation 2.4}
  \left\{
  \begin{aligned}
    u^*=u_a \qquad & \text{if} \;\; \alpha u^*-p>0, \\
    u^*\in \left[u_a,u_b\right] \qquad & \text{if} \;\; \alpha u^*-p=0, \\
    u^*=u_b \qquad & \text{if} \;\; \alpha u^*-p<0.
  \end{aligned}
  \right.
\end{equation}
By employing the weak minimum principle, we observe that the optimal control $u^*$ satisfies
\begin{equation*}
  \min_{v \in \left[u_a,u_b\right]} \left \{ \frac{\alpha}{2}v^2-pv \right \}=\frac{\alpha}{2}u^{*2}-pu^*.
\end{equation*}
From this, we see that
\begin{equation*}
  u^*=\frac{1}{\alpha}p
\end{equation*}
is true. 

\noindent Combining the above condition (\ref{equation 2.4}), we can get
\begin{equation}
  \label{equation 2.5}
  u^*=\Phi\left(\frac{p}{\alpha}\right)=\max\left\{u_a,\min\left\{u_b,\frac{p}{\alpha}\right\}\right\}. 
\end{equation}
Using (\ref{equation 2.3}) into (\ref{equation 2.2}), the optimality system (\ref{equation 2.2}) can be rewritten as:
\begin{equation}
  \label{equation 2.6}
  \begin{gathered}
    -\varDelta y+S(y)-\Phi\left(\frac{p}{\alpha}\right)=f  \quad \text{in}\;\Omega \quad \text{and} \quad y = 0  \quad \text{on}\;{\partial}{\Omega}, \\
    -\varDelta p+S'(y)p+y=y_d  \quad \text{in}\;\Omega \quad \text{and} \quad p = 0  \quad \text{on}\;{\partial}{\Omega}.
  \end{gathered}
\end{equation}

To solve the above problem (\ref{equation 2.6}) numerically, we use the second order five-point finite difference method to discretize (\ref{equation 2.6}). 
Then, discretizing $\Omega$ through an uniform Cartesian grid
\begin{equation*}
  \Omega_h=\left\{\left(x^i_1,x^j_2\right)=\left(ih,jh\right)| i=0,\dots,n+1; j=0,\dots,n+1.\right\}
\end{equation*}
with grid size $h=1/n$. Define $u_{ij}$, $y_{ij}$, $p_{ij}$, $f_{ij}$, $y_{dij}$ as an approximation of $u\left(x^i_1,x^j_2\right)$, 
$y\left(x^i_1,x^j_2\right)$, $p\left(x^i_1,x^j_2\right)$, $f\left(x^i_1,x^j_2\right)$, $y_d\left(x^i_1,x^j_2\right)$, respectively.
And let $u_h$, $y_h$, $p_h$, $f_h$, $y_{dh}$ be the corresponding lexicographic ordering (vectorization) of those approximations overall interior grid points.
Define $\Delta_h$ as a discretization of Laplacian operator $\Delta$, where $\Delta_h=(I\otimes J_h)+(J_h\otimes I)$ with $J_h=tridiag(1,-2,1)/h^2$ and I is identity matrix.
After discretization, we get the following solvable nonlinear equation
\begin{equation}
  \label{equation 2.7}
  \begin{gathered}
    F(y_h,p_h)=
        \begin{bmatrix}
            -\varDelta_h y_h+S(y_h)-\Phi(p/\alpha)-f_h \\
            -\varDelta_h p_h+S'(y_h)p_h+y_h-y_{dh}
        \end{bmatrix}=0,
  \end{gathered}
\end{equation}
where $S(\cdot)$, $S'(\cdot)$ and $\Phi (\cdot)$ are element-wisely defined and so is the multiplication $S'(y_h)p_h$.

\section{Inexact semismooth Newton-Krylov method with nonmonotonic line search}
\label{section 3}
The Rademacher theorem and definition of the generalized Jacobian are not held in function spaces. 
To address the nonsmoothness challenges in applying semismooth Newton methods to (\ref{equation 2.7}), we introduce the definition of slant differentiability. A merit function $Q(y_h,p_h)$ is constructed to quantify the solution quality of (\ref{equation 2.7}), which subsequently guides a nonmonotone line search strategy to ensure global convergence. For convenience, we denote the discrete $L^2$-norm as $\Vert \cdot\Vert $ in next.

\begin{Definition}
  \label{Definition 3.1}
    \textnormal{\cite{MR1786137}} A function $F:D \subset X \to Y$ is said to be slantly differentiable in an open domain $D_0 \subset D$
    if there exists a mapping $G:D\to L(X,Y)$ such that 
    \begin{equation}
      \label{equation 3.1}
      \lim_{h\to 0} \frac{\Vert  F(v+h)-F(v)-G(v+h)h\Vert }{\Vert  h \Vert }=0
    \end{equation}
    for every $v\in X$. Then $G$ is called a slanting function for $F$ in $D_0$.
\end{Definition}

\begin{Lemma}
    \label{Lemma 3.1} \textnormal{\cite{MR1786137}}An operator $F$ : $X\to Y$ is slantly differentiable at $x$ if and only if $F$ is Lipschitz continuous at $x$
\end{Lemma}
According to \cite{MR3461397}, it follows that $F(y_h,p_h)$ is slantly differentiable and $F(y_h,p_h)$ in (\ref{equation 2.7}) has a slanting function
\begin{equation}
  \label{equation 3.2}
  \begin{gathered}
    G_h(y_h,p_h) = 
    \begin{bmatrix}
      -\varDelta_h + \Psi(S'(y_h)) & -\frac{1}{\alpha}\Psi\left({\partial\Phi \left(\frac{1}{\alpha}p_h\right)}\right) \\
      I+\Psi(S''(y_h)p_h) & -\varDelta_h+\Psi(S'(y_h))
    \end{bmatrix}
  \end{gathered}
\end{equation}
where $\Psi(\cdot)$ denotes a diagonal matrix with the input vector as the diagonal elements. 
$\partial \Phi $ is chosen as the following form
\[ \partial \Phi(v) = \left\{
\begin{array}{rl}
  1 & \quad \text{if} \; u_a<v<u_b, \\
  0 & \quad \text{otherwise}.
\end{array} \right. \]

For the semismooth Newton (SSN) method to achieve optimal convergence performance, the slanting function must satisfy uniform boundness and invertibility conditions within open neighborhoods containing the optimal control-state pair. We formulate the following assumptions regarding the nonlinear terms in (\ref{equation 1.2}).

\begin{Assumption}
  \label{Assumption 3.1}
  Suppose the nonlinear functional $S(y)$ satisfies
  \begin{itemize}
    \item $S(y)\in C^3$ and $S'(y)$ is non-negative.
    \item $S(y)p+1\geq \kappa_1$ for some $\kappa_1>0$ in some neighborhood of the optimal $y^*$ and $p^*$.
  \end{itemize}
\end{Assumption}

\begin{Lemma}
  \label{Lemma 3.2}
  \textnormal{\cite{MR3461397}} Under the Assumption \ref{Assumption 3.1}, $G(y_h,p_h)$ is inverse and $\Vert  G(y_h,p_h)\Vert $ is uniformly bounded, there exists a constant $M>0$ such that
  \begin{equation}
    \label{equation 3.3}
    \Vert  G_h(y_h,p_h)^{-1} \Vert  \leq M,\Vert  G_h(y_h,p_h)\Vert  \leq M, \;\text{for all}\; h>0.
  \end{equation}
\end{Lemma}

\begin{Assumption}
    \label{Assumption 3.2} There exists a positive constant $L$ and $\delta_1$ such that for $\Vert  z_k - z^*\Vert <\delta_1$ and all sufficiently small $d$, we have
    \begin{equation}
        \Vert  G(z_k+d) - G(z_k) \Vert  \leq \frac{L\Vert  d\Vert }{M},
    \end{equation}
    where $M$ is mentioned in (\ref{equation 3.3}) and $z^*$ is an optimal solution.
\end{Assumption}

It is worth mentioning that the property discussed can be considered as the locally Lipschitz continuity of $G(y_h,p_h)$, which has been cited in several papers, such as \cite{MR1273766}. However, the function $F$ in this paper is not continuously differentiable as in \cite{MR1273766}, hence we treat this property as an assumption. Nevertheless, this assumption is valid when the grid discretization is sufficiently fine.

For ease of notation, we define 
$z_k:=[y_k;p_k]$, $F(z_k):=F(y_k,p_k)$ and $G(z_k):=G(y_k,p_k)$. We set the merit function $Q(z_k) = \frac{1}{2}\Vert  F(z_k)\Vert  ^2$. Next, we will prove that $Q(z_k)$ is continuously differentiable with gradient $\nabla Q(z_k) = G^T(z_k)F(z_k)$, where $G(z_k)$ is defined in (\ref{equation 3.2}). 

\begin{Lemma}
    \label{Lemma 3.3} Assume Assumption \ref{Assumption 3.1} and \ref{Assumption 3.2} are hold. For the function $F(z_k)$ given by (\ref{equation 2.7}), the merit function $Q(z_k)=\frac{1}{2}\Vert  F(z_k)\Vert  ^2$ is  continuously differentiable with gradient $\nabla Q(z_k) = G(z_k)^TF(z_k)$, where $G(z_k)$ is defined in (\ref{equation 3.2}). 
\end{Lemma}
\noindent {\it Proof. } The nonsmooth term of the function $F$ is denoted by $\Phi(\frac{p_h}{\alpha})$. According to the definition provided in equation (2.3), $F$ is continuously differentiable when $u_a<\frac{p_h}{\alpha}<u_b$, and in this case, $G(z_k)$ is equivalent to $\nabla F(z_k)$, hence 
\begin{equation}
    \label{equation 3.5} \nabla Q(z_k) = (\nabla F(z_k))^T F(z_k) = G(z_k)^TF(z_k).
\end{equation}
$F$ is not differentiable under the traditional gradient when $\frac{p_h}{\alpha}\leq u_a$ or $\frac{p_h}{\alpha}\geq u_b$. To address this, we utilize $G(z_k)$ as a substitute for $\nabla F(x)$. Let $F(z_k)=(F_1,...,F_{2(n-1)^2})$ and $G(z_k) = (G_1,...,G_{2(n-1)^2})^T$, according to the Lemma \ref{Lemma 3.1}, the component function $F_i$ of $F(z_k)$ is also Lipschitz continuous. 

When $u_a<\frac{p_h}{\alpha}<u_b$, we have
\begin{equation}
\label{equation 3.6}
    \nabla Q(z_k) = \nabla(\frac{1}{2}\Vert  F(z_k)\Vert  ^2)=\sum^{2(n-1)^2}_{i=1} F_i^2 G_i = G^T(z_k)F(z_k)
\end{equation}
To prove the continuity of $\nabla Q(z_k)$ at undifferentiable point, we have to prove $\Vert \nabla Q(z_k)\Vert \to \Vert \nabla Q(z_k +d)\Vert $ as $d\to 0$. Based on (\ref{equation 3.1}), we can deduce that
\begin{equation*}
    \Vert  F(z_k +d)-F(z_k)-G(z_k+d)d\Vert  = o(\Vert  d\Vert ).
\end{equation*}
Combing with Lemma \ref{Lemma 3.2} we get 
\begin{align*}
    \Vert  G(z_k+d)^T F(z_k+d)-G(z_k+d)^TF(z_k)\Vert  &\leq \Vert  G(z_k+d)^T(F(z_k+d)-F(z_k))\Vert \\
    &\leq \Vert   G(z_k+d)^T(G(z_k+d)d+o(\Vert d \Vert))\Vert\\
    &\leq  M^2\Vert  d\Vert ^2+o(\Vert  d \Vert ) 
\end{align*}
Thus
\begin{equation}
    \label{equation 3.7} \lim_{d\to 0} \Vert  G(z_k+d)^TF(z_k+d)\Vert  = \lim_{d\to 0} \Vert  G(z_k+d)^TF(z_k)\Vert .
\end{equation}
Moreover,
\begin{equation}
\label{equation 3.8}
    \Vert  G(z_k+d)^TF(z_k)- G(z_k)^TF(z_k)\Vert = \Vert ( G(z_k+d)^T- G(z_k)^T)F(z_k)\Vert .
\end{equation}
From Lemma \ref{Lemma 3.2} and Assumption \ref{Assumption 3.2}, there exists a $\delta_1$ and $L$ such that
\begin{equation}
\label{equation 3.9}
    \Vert ( G(z_k+d)^T- G(z_k)^T)F(z_k)\Vert  \leq \frac{L\Vert  F(z_k)\Vert  }{M}\Vert  d\Vert .
\end{equation}
Combing (\ref{equation 3.7}), (\ref{equation 3.8}) and (\ref{equation 3.9}), we have
\begin{equation}
    \lim _{d\to0}\Vert  G(z_k)^TF(z_k)\Vert  = \lim _{d\to0}\Vert  G(z_k+d)^TF(z_k)\Vert  = \lim_{d\to 0} \Vert  G(z_k+d)^TF(z_k+d)\Vert 
\end{equation}

In conclusion, whether $u_a<\frac{p_h}{\alpha}<u_b$, $\frac{p_h}{\alpha}\leq u_a$ or $\frac{p_h}{\alpha}\geq u_b$, we have $\Vert \nabla Q(z_k)\Vert \to \Vert \nabla Q(z_k +d)\Vert $ as $d\to 0$. The proof is complete.

There are numerous strategies to achieve global convergence of the Newton method, such as trust region methods, line search, and backtracking, among others. In this paper, we adopt the nonmonotone line search strategy, which can be found in \cite{MR849278}, that is
\begin{equation*}
    Q(z_k + \delta_kd_k)\leq \underset{0\leq j \leq k }{max}Q(z_j)+c_1\delta_k\nabla Q(z_k)^Td_k.
\end{equation*}

In each iteration of the Newton method, obtaining  $d_k$ requires solving the equation $G(z_k)d_k=-F(z_k)$.  However, precisely solving $d_k$ is computationally expensive and involves inverting large-scale matrices. To avoid large-scale matrix inversion and reduce computational costs, we use the Generalized Minimum Residual Method(GMRES) to obtain $d_k$ inexactly, satisfying:
\begin{equation*}
    \Vert  F(z_k)+G(z_k)d_k\Vert \leq \eta_k\Vert  F(z_k)\Vert ,
\end{equation*}
where $\eta_k\in [0,\eta_{max}],\eta_{max}\in [0,1)$.

For the function $F$ defined by (\ref{equation 2.7}), we have devised the Inexact semismooth Newton-GMRES method with nonmonotone line search (ISSNG-L) as follows,  the stopping criterion $\tau_k$  is specified in the experimental section:
\begin{algorithm}[H]
  \caption{ISSNG-L}
  \label{GC-SSN}
  \begin{algorithmic}[1]
    \Require Given constant $c_1\in (0,1),\delta_0\in(0,1],\eta_{max}\in [0,1)$ and $0<\theta_{min}<\theta_{max}<1$, allowable error $\epsilon>0$, initialize $z_0$, set $k=0$. 
    \Ensure Final numerical solution $\bar{z}=z_{k+1}$. 
    \State Compute corresponding $G(z_k)$ and $F(z_k)$ as in (\ref{equation 3.2}) and (\ref{equation 2.7}). Find some $\eta_k\in [0,\eta_{max}]$ and use the GMRES method to solve the Newton step equation $G(z_k)d_k=-F(z_k)$ to obtain $d_k$ that satisfies
    \begin{equation}
    \label{equation 3.11}
        \Vert  F(z_k)+G(z_k)d_k\Vert \leq \eta_k\Vert  F(z_k)\Vert 
    \end{equation}
    \State Choose $\theta\in[\theta_{min},\theta_{max}]$
    \While {$ Q(z_k + \delta_kd_k)>\underset{0\leq j \leq k }{max}Q(z_j)+c_1\delta_k\nabla Q(z_k)^Td_k$}
      \State Set $\delta_k=\theta \delta_k$
    \EndWhile
    \State Let $z_{k+1}=z_k+\delta_k d_k$ and $k=k+1$.
    \If {$\tau_k<Tol$}
      \State {Output $z_{k+1}$.}
    \Else
      \State Return step 1.
    \EndIf
  \end{algorithmic}
\end{algorithm}

\begin{Remark}
\end{Remark}

The adoption of the GMRES method for inexactly solving the Newton equation $G(z_k)d_k=-F(z_k)$ is primarily motivated by its superior efficiency in both storage and computational cost compared to direct solution methods, especially for large-scale problems arising\cite{MR2278446,MR1273766,MR1860971}:

(Storage Complexity)With the problem size $N=2(n-1)^2$ in this paper, the direct approach (based on LU factorization) explicitly stores the entire Jacobian matrix in memory, leading to a storage complexity of $O(N^2)=O(n^4)$. The GMRES iteration adopts a matrix-free strategy that avoids storing the full Jacobian $G(z_k)$, it keeps only a Krylov basis of dimension $m$ (typically $m<<  N$), yielding a storage complexity of $O(mN)=O(mn^2)=O(n^2)$. 

(Computational Complexity)With the problem size $N=2(n-1)^2$ in this paper, The direct approach (based on LU factorization) incurs a computational complexity of $O(N^{\frac{3}{2}})=O(n^3)$. By contrast, the overall computational complexity of GMRES iteration is $O(mN+m^2N)=O(n^2)$.

The preceding complexity analysis demonstrates that GMRES iteration achieves significantly higher computational efficiency and lower cost than the direct approach, this observation is corroborated by the numerical comparisons reported in Example 1 of Section \ref{section 5}.

\section{Convergence analysis}
\label{section 4}
In this section, we first demonstrate that $d_k$, satisfying (\ref{equation 3.11}), constitutes a descent direction for $Q(z_k)$. Subsequently, we establish the global convergence of Algorithm \ref{GC-SSN}. Finally, We have demonstrated that Algorithm \ref{GC-SSN} achieves superlinear convergence when $\eta_k\to0$ in (\ref{equation 3.11}).

Given $z_0$, we denote a sequence ${z_k}$ generated by Algorithm 1 and the level set of $\Vert  F(z_k)\Vert $ as
\begin{equation}
\label{equation 4.1}
    L(z_0)=\{z_k|\Vert  F(z_k)\Vert  \leq \Vert  F(z_0)\Vert  \}.
\end{equation}
\begin{Assumption}
    \label{Assumption 4.1} The sequence $\{z_k\}$ generated by Algorithm 1 is contained in a compact set $L(z_o)$.
\end{Assumption}

\begin{Lemma}
    \label{Lemma 4.1} Assume that there exists a direction $d_k$ that satisfies (\ref{equation 3.11}) under the condition $\Vert  G(z_k)^TF(z_k)\Vert =0$, then $\Vert  F(z_k)\Vert =0$.
\end{Lemma}
\noindent {\it Proof. } According to the inequality (\ref{equation 3.11}), we have
\begin{equation}
    \label{equation 4.2} \Vert  F(z_k)+G(z_k)d_k\Vert ^2 = \Vert  F(z_k)\Vert ^2 + 2[(G(z_k)^TF(z_k)]^Td_k+\Vert  G(z_k)d_k\Vert ^2\leq \eta_k^2\Vert  F(z_k)\Vert ^2.
\end{equation}
Since $\Vert  G(z_k)^TF(z_k)\Vert =0$ and $\eta_k\in (0,1)$, we have
\begin{equation*}
    \Vert  G(z_k)d_k\Vert ^2\leq -(1-\eta_k^2)\Vert  F(z_k)\Vert ^2\leq 0,
\end{equation*}
Hence, $\Vert  F(z_k)\Vert =0$. The proof is complete.


\begin{Lemma}
    \label{Lemma 4.2} Assume that $d_k$ satisfies (\ref{equation 3.11}), then $d_k$ is descent direction for $Q(z_k)$ at the current approximation $z_k$ such that
    \begin{align}
        \label{equation 4.3}&-\nabla Q(z_k)^Td_k\geq (1-\eta_k)\Vert  F(z_k)\Vert ^2>0\\
   \label{equation 4.4}  &\frac{|\nabla Q(z_k)^Td_k|}{\Vert  d_k\Vert } \geq \frac{1-\eta_k}{M^2(1+\eta_k)}\Vert  \nabla Q(z_k)\Vert \geq0
    \end{align}
where $M$ is given in Lemma \ref{Lemma 3.2} and $\eta_k\in[0,1)$.
\end{Lemma}
\noindent{\it Proof. } For the linear system $G(z_k)d_k=-F(z_k)$, We set $r_k=F(z_k)+G(z_k)d_k$ as the residual, and according to the Lemma \ref{Lemma 3.3} and (\ref{equation 3.11}), we have
\begin{align}
    \nabla Q(z_k)d_k= (G(z_k)^TF(z_k))^Td_k=F(z_k)^TG(z_k)dk&=F(z_k)^T(r_k-F(z_k))\nonumber\\
    &\leq \Vert  F(z_k)\Vert \cdot \Vert  r_k\Vert  -\Vert  F(z_k)\Vert ^2\nonumber\\
    &\leq \eta_k\Vert  F(z_k)\Vert  \cdot \Vert  F(z_k)\Vert  - \Vert  F(z_k)\Vert \nonumber\\
  \label{equation 4.5}   &=-(1-\eta_k)\Vert  F(z_k)\Vert  ^2,
\end{align}
which proves (\ref{equation 4.3}) holds true.

It is clearly that if $G(z_k)^TF(z_k)\not=0$, then $d_k\not=0$, it indicates that
\begin{align}
    \Vert  d_k\Vert  &=\Vert  G(z_k)^{-1}(r_k-F(z_k))\Vert \nonumber\\
    &\leq \Vert  G(z_k)^{-1}\Vert  (\Vert  r_k\Vert +\Vert  F(z_k)\Vert )\nonumber\\
    &\leq \Vert  G(z_k)^{-1}\Vert  (\eta_k\Vert  F(z_k) \Vert +\Vert  F(z_k)\Vert )\nonumber\\
   \label{equation 4.6} &=(1+\eta_k)\Vert  G(z_k)^{-1}\Vert  \Vert  F(z_k)\Vert .
\end{align}
Incorporating (\ref{equation 4.5}), (\ref{equation 4.6}) and the fact $\Vert  \nabla Q(z_k)\Vert  \leq \Vert  G(z_k)\Vert \Vert  F(z_k)\Vert $, we have
\begin{align}
    \frac{\Vert  \nabla Q(z_k)^Td_k\Vert }{\Vert  \nabla Q(z_k)\Vert  \Vert  d_k\Vert } &\geq \frac{(1-\eta_k)\Vert  F(z_k)\Vert ^2}{(1+\eta_k)\Vert  G(z_k)^{-1}\Vert  \Vert  F(z_k)\Vert \Vert  \nabla Q(z_k)\Vert }\nonumber\\
    &=\frac{(1-\eta_k)\Vert  F(z_k)\Vert }{(1+\eta_k)\Vert  G(z_k)^{-1}\Vert \Vert  \nabla Q(z_k)\Vert }\nonumber\\
    &\frac{(1-\eta_k)\Vert  F(z_k)\Vert }{(1+\eta_k)\Vert  G(z_k)^{-1}\Vert \Vert  G(z_k)\Vert \Vert  F(z_k)\Vert }\nonumber\\
    \label{equation 4.7}&=\frac{1-\eta_k}{(1+\eta_k)\Vert  G(z_k)^{-1}\Vert \Vert  G(z_k)\Vert }.
\end{align}
From Lemma \ref{Lemma 3.2}, we have $\Vert  G(z_k)\Vert \leq M, \Vert  G(z_k)^{-1}\Vert \leq M$, and combing it with (\ref{equation 4.7}), we obtain
\begin{equation}
    \label{equation 4.8} \frac{\Vert  \nabla Q(z_k)^Td_k\Vert }{\Vert  \nabla Q(z_k)\Vert  \Vert  d_k\Vert } \geq \frac{1-\eta_k}{(1+\eta_k)M^2}.
\end{equation}
Multiplying both sides of (\ref{equation 4.8}) by $\nabla Q(z_k)$, we get
\begin{equation}
    \label{equation 4.9} \frac{|\nabla Q(z_k)^Td_k|}{\Vert  d_k\Vert } \geq \frac{1-\eta_k}{M^2(1+\eta_k)}\Vert  \nabla Q(z_k)\Vert \ >0,
\end{equation}
which proves (\ref{equation 4.4}) holds true. The proof is complete.

\begin{Lemma}
    \label{Lemma 4.3} 
    \textnormal Suppose $d_k$ is a direction that satisfies (\ref{equation 3.11}), $\nabla Q(z_k)\not=0$, $\eta_k\in[0,\eta_{max}]$ and $\eta_{max}\in[0,1)$. Then, there must exist a $\delta_k$ such that after a finite number of backtracking steps, the following condition is satisfied:
    \begin{equation}
       \label{equation 4.10}Q(z_k + \delta_kd_k)\leq \underset{0\leq j \leq k }{max}Q(z_j)+c_1\delta_k\nabla Q(z_k)^Td_k.
    \end{equation}
\end{Lemma}
\noindent {\it Proof. } From Lemma (\ref{Lemma 3.3}), we know $\nabla Q(z_k)$ is continuously differentiable, we can apply the Lagrange Mean Value Theorem to obtain:
\begin{equation}
\label{equation 4.11}
    \nabla Q(z_k+v_k\delta_kd_k)=\frac{Q(z_k+\delta_kd_k)-Q(z_k)}{\delta_kd_k},
\end{equation}
where $v_k\in[0,1]$. By rearranging (\ref{equation 4.11}),  we obtain:
\begin{align}
    &Q(z_k+\delta_kd_k)\nonumber\\
    =&Q(z_k)+c_1\delta_k\nabla Q(z_k)^Td_k-c_1\delta_k\nabla Q(z_k)^Td_k+\delta_k\nabla Q(z_k)^Td_k-\delta_k\nabla Q(z_k)^Td_k+\delta_k\nabla Q(z_k+v_k\delta_kd_k)^Td_k\nonumber\\
    =&Q(z_k)+c_1\delta_k\nabla Q(z_k)^Td_k+(1-c_1)\delta_k\nabla Q(z_k)^Td_k+\delta_k[\nabla Q(z_k+v_k\delta_kd_k)^Td_k-\nabla Q(z_k)^Td_k]\nonumber\\
    =&Q(z_k)+c_1\delta_k\nabla Q(z_k)^Td_k+\delta_k[\nabla Q(z_k+v_k\delta_kd_k)^Td_k-\nabla Q(z_k)^Td_k+(1-c_1)\nabla Q(z_k)^Td_k]\nonumber\\
    \label{equation 4.12}\leq &\underset{1\leq j\leq k}{max}Q(z_j)+c_1\delta_k\nabla Q(z_k)^Td_k+\delta_k[\nabla Q(z_k+v_k\delta_kd_k)^Td_k-\nabla Q(z_k)^Td_k+(1-c_1)\nabla Q(z_k)^Td_k]
\end{align}
From Lemma (\ref{Lemma 3.3}), we know $\nabla Q(z_k)$ is continuous. Since $\nabla Q(z_k)\not=0$, there exists $\mu$ and $\sigma$, for all $\Vert  z_k-z\Vert \leq \sigma$, we have $\nabla Q(z)\geq\mu$. In addition, there exists a sufficiently small stepsize $\delta_k$ such that
\begin{equation}
    \label{equation 4.13} \Vert  \nabla Q(z_k+v_k\delta_k)^T-\nabla Q(z_k)^T\Vert  \leq
    (1-c_1)\frac{1-\eta_{max}}{M^2(1+\eta_{max})}\mu,
\end{equation}
where $M$ is defined in Lemma \ref{Lemma 3.2}, $\eta_k\in[0,\eta_{max}]$ and $\eta_{max}\in[0,1)$. Multiplying both sides of the above inequality by $\Vert  d_k\Vert $ yields:
\begin{equation}
\label{equation 4.14}
    |[\nabla Q(z_k+v_k\delta_k)^Td_k-\nabla Q(z_k)^T]d_k|\leq\Vert  \nabla Q(z_k+v_k\delta_k)^T-\nabla Q(z_k)^T\Vert  \Vert  d_k\Vert  \leq \frac{(1-c_1)(1-\eta_{max})}{M^2(1+\eta_{max})}\mu \Vert  d_k\Vert .
\end{equation}
By (\ref{equation 4.3}) and (\ref{equation 4.4}), we have
\begin{equation}
    \label{equation 4.15} \nabla Q(z_k)^Td_k\leq-\frac{1-\eta_k}{M^2(1+\eta_k)}\Vert  \nabla Q(z_k)\Vert  \Vert  d_k\Vert \leq -\frac{1-\eta_{max}}{M^2(1+\eta_{max})}\mu\Vert  d_k\Vert 
\end{equation}
Based on (\ref{equation 4.12}), (\ref{equation 4.14}) and (\ref{equation 4.15}), we can conclude that the value within the brackets of the third term on the right-hand side of (\ref{equation 4.12}) will eventually become negative after a finite number of reductions. At this stage, $\delta_k$ becomes admissible, implying that there exists an $\delta_k$ satisfying the following inequality after a finite number of backtracking steps:
\begin{equation}
    Q(z_k + \delta_kd_k)\leq \underset{1\leq j\leq k}{max}Q(z_j)+c_1\delta_k\nabla Q(z_k)^Td_k.
\end{equation}
The proof is complete.

Before analyzing the global convergence and local superlinear convergence rate of Algorithm \ref{GC-SSN}, we first present a well-known theorem established in \cite{MR849278}, the parameter notations in this theorem are independent of those used in our work.

\begin{Lemma}
    \label{Lemma 4.4} \textnormal{\cite{MR849278}} Let $\{x_k\}$ be a sequence defined by 
    \begin{equation*}
        x_{k+1}=x_k+\alpha_kd_k, d\not=0.
    \end{equation*}
Let $\alpha>0$, $\sigma\in(0,1)$, $\gamma\in(0,1)$ and let $M$ be a nonnegative integer. Assume that:
\begin{enumerate}[label=(\roman*)]
    \item the level set $\Omega_0=\{x:f(x)\leq f(x_0)\}$ is compact;
    \item there exists positive numbers $c_1$, $c_2$ such that:
    \begin{align}
        \nabla Q(x_k)^Td_k&\leq -c_1\Vert  \nabla Q(x_k)\Vert  ^2\\
        \Vert  d_k\Vert  &\leq c_2\Vert  \nabla Q(x_k)\Vert ;
    \end{align}
    \item  $\alpha_k=\sigma^{h_k}$, where $h_k$ is the first nonnegative integer $h$ for which:
    \begin{equation}
        Q(x_k+\lambda_kd_k)\leq \underset{0\leq j\leq m(k)}{\max} [Q(x_{k-j})]+\gamma\omega_k\alpha\nabla Q(x_k)^Td_k
    \end{equation}
    where $m(0)=0$ and $0\leq m(k)\leq\min[m(k-1)+1,M]$, $k\geq1$.
\end{enumerate}
Then:
\begin{enumerate}[label=(\alph*)]
    \item the sequence \(\{x_k\}\) remains in \(\Omega_0\) and every limit point \(\bar{x}\) satisfies \(g(\bar{x}) = 0\);
    \item no limit point of \(\{x_k\}\) is a local maximum of \(f\);
    \item if the number of the stationary points of \(f\) in \(\Omega_0\) is finite, the sequence \(\{x_k\}\) converges.
\end{enumerate}
\end{Lemma}

We can obtain the following conclusions from Lemma \ref{Lemma 4.4}.

\begin{Theorem}
    \label{Theorem 4.1} Consider a sequence $\{z_k\}$ produced by Algorithm \ref{GC-SSN}. Set $d_k$ fulfill (\ref{equation 3.11}). Assume the existence of constant $\eta_{max}\in[0,1)$ such that $\eta_k\leq\eta_{max}$ in (\ref{equation 3.11}), Assumption \ref{Assumption 3.1}, Assumption \ref{Assumption 3.2} and Assumption \ref{Assumption 4.1} hold. Suppose $\nabla Q(z_k)\not=0$. 
    Then:
    \begin{enumerate}[label=(\alph*)]
    \item the sequence  \(\{z_k\}\) stays within the level set $L(z_0)$ and every limit point $z^*$ satisfies $Q(z_k)=0$;  
    \item  if the number of the stationary points of \(Q(z_k)\) in \(\Omega_0\) is finite, the sequence \(\{z_k\}\) converges.
    \end{enumerate}
\end{Theorem}
\noindent{\it Proof. } The proof will proceed by verifying each of the three cases in Lemma \ref{Lemma 4.4} in turn:
\begin{enumerate}[label=(\roman*)]
    \item From Assumption \ref{Assumption 4.1}, we deduce that the sequence $\{z_k\}$ generated by Algorithm \ref{GC-SSN} is contained in a compact set $L(z_o)$, which confirms the first condition;
    \item By the $\nabla Q(z_k)=G(z_k)^TF(z_k)$ in Lemma \ref{Lemma 3.3} and Lemma \ref{Lemma 3.2}, we have
    \begin{equation}
        \label{equation 4.20} \Vert  F(z_k)\Vert ^2 = \frac{\Vert  \nabla Q(z_k)\Vert  ^2}{\Vert  (G(z_k)^T)^{-1}\Vert  ^2}= \frac{\Vert  \nabla Q(z_k)\Vert  ^2}{\Vert  (G(z_k)^{-1})^T\Vert  ^2} =\frac{\Vert  \nabla Q(z_k)\Vert  ^2}{\Vert  (G(z_k)^{-1}\Vert  ^2}\geq \frac{\Vert  \nabla Q(z_k)\Vert  ^2}{M^2}.
    \end{equation}
    Combining (\ref{equation 4.20}) and (\ref{equation 4.5}), we obtain
    \begin{equation}
        \label{equation 4.21} \nabla Q(z_k)^Td_k\leq -(1-\eta_k)\Vert  F(z_k)\Vert  ^2 \leq -(1-\eta_{max})\Vert  F(z_k)\Vert  ^2 \leq -\frac{(1-\eta_{max})}{M^2}\Vert  \nabla Q(z_k)\Vert^2 .
    \end{equation}
    From (\ref{equation 4.6}), we get
    \begin{equation}
        \label{equation 4.22} \Vert  d_k\Vert  \leq (1+\eta_k)\Vert  G(z_k)^{-1}\Vert  \Vert  F(z_k)\Vert \leq(1+\eta_{max})M \Vert  F(z_k)\Vert .
    \end{equation}
    (\ref{equation 4.21}) and (\ref{equation 4.22}) together constitute the second condition, thereby the second condition holds.
    \item The third condition is satisfied with the sequence of backtracking parameters $\delta_k=1,\omega,\omega^2,...$ in
    \begin{equation}
        \label{equation 4.23}  Q(z_k + \delta_kd_k)\leq\underset{1\leq j\leq k}{max}Q(z_j)+c_1\delta_k\nabla Q(z_k)^Td_k.
    \end{equation}
\end{enumerate}

In summary, since all three conditions of Lemma \ref{Lemma 4.4} have been verified, the conclusion of the Theorem holds.

\begin{Theorem}
    \label{Theorem 4.2} Consider the sequence $\{z_k\}$ generated by Algorithm \ref{GC-SSN} under the conditions specified in Theorem \ref{Theorem 4.1}. If $c_1\leq \frac{1}{2}$ and $\eta_k\to0$, then the stepsize $\delta_k=1$ is accepted in (\ref{equation 4.10}) for sufficiently large $k$.
\end{Theorem}
\noindent{\it Proof. } According to the $Q(z_k)=\frac{1}{2}\Vert  F(z_k)\Vert ^2$, we have
\begin{align}
    Q(z_k+d_k)-Q(z_k)-\nabla Q(z_k)^Td_k &=\frac{1}{2}\Vert  F(z_k+d_k)\Vert  - \frac{1}{2}\Vert  F(z_k)\Vert  -\nabla Q(z_k)^Td_k\nonumber\\
    &=\frac{1}{2}\Vert  F(z_k)+G(z_k)d_k+o(\Vert  d_k\Vert )\Vert  ^2- \frac{1}{2}\Vert  F(z_k)\Vert^2  -\nabla Q(z_k)^Td_k\nonumber\\
    \label{equation 4.24}&=\frac{1}{2}\Vert  G(z_k)d_k\Vert  ^2+o(\Vert  d_k\Vert ^2),
\end{align}
which implies
\begin{equation}
     \label{equation 4.25}  Q(z_k+d_k)=Q(z_k)+c_1\nabla Q(z_k)^Td_k +\frac{1}{2}(\nabla Q(z_k)^Td_k+\Vert  G(z_k)d_k\Vert ^2)+(\frac{1}{2}-c_1)\nabla Q(z_k)^Td_k+o(\Vert  d_k\Vert ^2).
\end{equation}
From (\ref{equation 4.3}), we get
\begin{equation}
\label{equation 4.26}
    \nabla Q(z_k)^Td_k\leq -(1-\eta_k)\Vert  F(z_k)\Vert ^2.
\end{equation}
By $G(z_k)d_k=r_k-F(z_k)$ and (\ref{equation 3.11}), we obtain
\begin{equation}
    \label{equation 4.27} \Vert  G(z_k)d_k\Vert  \leq \Vert  r_k\Vert  +\Vert  F(z_k)\Vert \leq (1+\eta_k)\Vert  F(z_k)\Vert .
\end{equation}
As (\ref{equation 4.22}) means
\begin{equation}
\label{equation 4.28}
    \Vert  d_k\Vert  \leq (1+\eta_k)\Vert  G(z_k)^{-1}\Vert  \Vert  F(z_k)\Vert \leq \Vert  d_k\Vert  \leq (1+\eta_k)M \Vert  F(z_k)\Vert ,
\end{equation}
we have $o(\Vert  d_k\Vert ^2)=o(\Vert  F(z_k)\Vert ^2)$.
Combining (\ref{equation 4.26}), (\ref{equation 4.27}) and (\ref{equation 4.28}), (\ref{equation 4.25}) can be rewritten as
\begin{align}
    Q(z_k+d_k)=&Q(z_k)+c_1\nabla Q(z_k)^Td_k +\frac{1}{2}(\nabla Q(z_k)^Td_k+\Vert  G(z_k)d_k\Vert ^2)+(\frac{1}{2}-c_1)\nabla Q(z_k)^Td_k+o(\Vert  d_k\Vert ^2)\nonumber\\
    \leq& Q(z_k)+c_1\nabla Q(z_k)^Td_k+\frac{1}{2}[ -(1-\eta_k)\Vert  F(z_k)\Vert ^2+(1+\eta_k)^2\Vert  F(z_k)\Vert ^2]\nonumber\\
    &-(\frac{1}{2}-c_1) (1-\eta_k)\Vert  F(z_k)\Vert ^2+o(\Vert  d_k\Vert ^2)\nonumber\\
    \leq &Q(z_k)+c_1\nabla Q(z_k)^Td_k-[(\frac{1}{2}-c_1) (1-\eta_k)+\frac{1}{2}(1-\eta_k)+(1+\eta_k)^2]\Vert  F(z_k)\Vert ^2+o(\Vert  d_k\Vert ^2)\nonumber\\
    \label{equation 4.29}\leq &Q(z_k)+c_1\nabla Q(z_k)^Td_k.
\end{align}
For sufficiently large $k$, the inequality holds. This is derived from the fact that the third term in the brackets on the right-hand side of (\ref{equation 4.29}) become negative, as shown by
\begin{equation*}
    (\frac{1}{2}-c_1) (1-\eta_k)+\frac{1}{2}(1-\eta_k)+(1+\eta_k)^2\to\frac{1}{2}-c_1
\end{equation*}
when $\eta_k\to0$. (\ref{equation 4.29}) implies the admissibility of $\delta_k=1$  for sufficiently large k in (\ref{equation 4.10}):
\begin{equation}
    \label{equation 4.30} z_{k+1}= z_k+d_k.
\end{equation}
The proof is complete.

In Theorem \ref{Theorem 4.2},  for sufficiently large $k$, when $z_k$ is sufficiently close to the optimal solution $z^*$, the Inexact Newton method will proceed with $\delta_k=1$ at each iteration. More precisely, there exists a neighborhood of $z^*$ where this holds. Since this characterization is not utilized in Theorem \ref{Theorem 4.2}, we omit it there but will incorporate it in the subsequent theorem.

\begin{Theorem}
    \label{Theorem 4.3}  (local superlinear convergence) Consider the sequence $\{z_k\}$ generated by Algorithm \ref{GC-SSN} under the conditions specified in Theorem \ref{Theorem 4.2}. Suppose $z^*$ is a limit point of $\{z_k\}$. There exists a positive constant $\rho$ such that for all $\Vert  z_k-z^*\Vert <\rho $, the sequence $\{z_k\}$ superlinear converges to $z^*$.
\end{Theorem}
\noindent{\it Proof. } When $z^*$ is a limit point of $\{z_k\}$, Theorem \ref{Theorem 4.2} gives
\begin{equation}
    \label{equation 4.31} F(z^*)=0.
\end{equation}
By (\ref{Definition 3.1}), there exists a positive constant $\rho$ such that for all $\Vert  z_k-z^*\Vert <\rho $ and $\forall \epsilon>0$, we have
\begin{equation}
    \label{equation 4.32} \Vert  F(z_k)-F(z^*)-G(z_k)(z_k-z^*)\Vert \leq \frac{\varepsilon}{M}\Vert  z_k-z^*\Vert ,
\end{equation}
where $M$ is given in Lemma \ref{Lemma 3.2}.
Furthermore, we have
\begin{align}
    \Vert z_k+d_k-z^*\Vert &=\Vert G^{-1}(z_k)[G(z_k)(z_{k}-z^*)+G(z_k)s_k]\Vert \nonumber\\
  &\leq \Vert G^{-1}(z_k)\Vert (\Vert G(z_k)(z_{k}-z^*)+F(z_k)-F(z_k)+G(z_k)d_k\Vert ) \nonumber\\
  &\leq  \Vert G^{-1}(z_k)\Vert[\Vert F(z^*)+G(z_k)(z_{k}-z^*)-F(z_k)\Vert +\Vert r_k\Vert]\nonumber\\
  &\leq \Vert G^{-1}(z_k)\Vert[\Vert F(z^*)+G(z_k)(z_{k}-z^*)-F(z_k)\Vert +\eta_k\Vert F(z_k)\Vert ]\nonumber\\
  &\leq \Vert G^{-1}(z_k)\Vert[\frac{\varepsilon}{M}\Vert z_k-z^*\Vert +\eta_k(\Vert G(z_k)(z_k-z^*)\Vert +\frac{\epsilon}{M}\Vert z_k-z^*\Vert )]\nonumber\\
  \label{equation 4.33}&\leq \epsilon \Vert z_k-z^*\Vert + \eta_k M^2\Vert z_k-z^*\Vert +\eta_k\epsilon\Vert z_k-z^*\Vert.
\end{align}
For $\forall \epsilon>0$, let $\eta_k\leq\frac{\epsilon}{M^2}<1$. Combining Lemma \ref{Lemma 3.2}, (\ref{equation 4.30}) and (\ref{equation 4.33}), we obtain
\begin{equation}
    \label{equation 4.34}\Vert z_{k+1}-z^*\Vert \leq  \epsilon \Vert z_k-z^*\Vert+ \epsilon^2\Vert z_k-z^*\Vert +\frac{\epsilon^2}{M^2}\Vert z_k-z^*\Vert  \leq (\epsilon+\epsilon^2+\frac{\epsilon^2}{M^2})\Vert z_k-z^*\Vert.
\end{equation}
$\eta_k=\frac{\epsilon}{M^2}\to0$ implies that $(\epsilon+\epsilon^2+\frac{\epsilon^2}{M^2})\to0$, so it is obviously that $z_{k+1}\in B_\rho(z^*)$. Hence ,from (\ref{equation 4.35}) we have
\begin{equation}
    \label{equation 4.35}\underset{k\to\infty}{lim} \frac{\Vert z_{k+1}-z^*\Vert}{\Vert z_k-z^*\Vert } =0,
\end{equation}
this theorem is proved.

\begin{Theorem}
    \label{Theorem 4.4} (global convergence)  Consider the sequence $\{z_k\}$ generated by Algorithm \ref{GC-SSN} under the conditions specified in Theorem \ref{Theorem 4.3}. Suppose $z^*$ is a limit point of $\{z_k\}$. Then the sequence $\{z_k\}$ converges to $z^*$ and $\Vert F(z^*)\Vert =0$.
\end{Theorem}
\noindent{\it Proof. } According to Theorem \ref{Theorem 4.3}, we know that the sequence $\{z_k\}$ generated by Algorithm \ref{GC-SSN} converges to $z^*$. From Theorem \ref{Theorem 4.1}, we can get $\Vert F(z^*)\Vert =0$. The proof is thus established.

\section{Numerical experiments}
\label{section 5}
In this section, we numerically illustrate the properties of the proposed the inexact semismooth Newton-GMRES method with nonmonotonic line search (ISSNG-L) and the inexact semismooth Newton-GMRES methods (ISSNG) through two numerical examples. We not only compare ISSNG-L with ISSNG to highlight the advantages of the nonmonotonic line search in ISSNG-L but also compare our method with similar method from other paper, demonstrating the superiority of our methods. All numerical results were obtained on a desktop computer equipped with an AMD Ryzen Threadripper PRO 5975WX processor featuring 32 cores, along with 256GB of RAM, utilizing MATLAB. We set $r_y$ and $r_p$ be the residual of state equation and costate equation:
\begin{equation*}
    r_y = -\varDelta_h y_h+S(y_h)-\Phi(p/\alpha)-f_h , r_p = -\varDelta_h p_h+S'(y_h)p_h+y_h-y_{dh}.
\end{equation*}
The stopping criterion is established as follows:
\begin{equation}
    \tau_k=\frac{\Vert r^k_y\Vert + \Vert r^k_p\Vert}{max(1,\Vert r^0_y\Vert + \Vert r^0_y\Vert)} \leq 10^{-8},
\end{equation}
where $r^k_y$ and $r^k_p$ are the residuals at $k$th iteration. In the following examples, we choose
\begin{equation}
    \eta_k = \gamma (\frac{\Vert F(z_k)\Vert }{\underset{1\leq j \leq k-1}{max}\Vert F(z_j)\Vert })^{a_1}, k=1,2...
\end{equation}
where $\gamma\in[0,1]$, $a_1\in(1,2]$ and $\eta_0\in[0,1)$.

\subsection{Example 1}

\cite{MR4447697} Let$S(y)= y^3$, $u_a=-\infty$, $u_b=\infty$, $\alpha=10^{-3}$, 
\begin{align*}
&z(x_1,x_2)= sin(\pi x_1)sin(\pi x_2), \\
&f=2\pi^2 z(x_1,x_2)+\frac{\pi^2}{10^3}z(x_1,x_2)-z(x_1,x_2)e^{\pi x_1}, \\
&y_d=z(x_1,x_2)+(\pi^2 z(x_1,x_2)e^{\pi x_1})/10^3 - (2\pi^2cos(\pi x_1)sin(\pi x_2)e^{\pi x_1})/10^3 +(3z^3(x_1,x_2)e^{\pi x_1})/10^3.
\end{align*}
The example possesses a unique optimal control $u^*=z(x_1,x_2)e^{\pi x_1}$. In the computations for this example, the initial guess is set to matrices filled entirely with zeros.

\begin{figure}[ht]
    \centering
    \begin{subfigure}{0.4\linewidth}
        \centering
        \includegraphics[width=\textwidth]{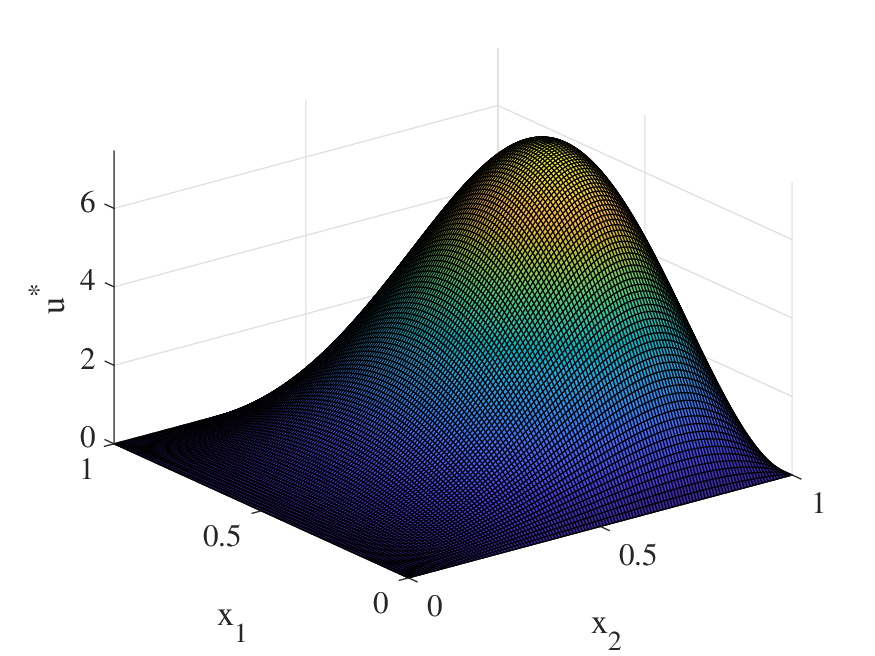}
        \caption{Exact control $u^*$}
        \label{fig:exact-control}
    \end{subfigure}
    \hspace{5pt} 
    \begin{subfigure}{0.4\linewidth}
        \centering
        \includegraphics[width=\textwidth]{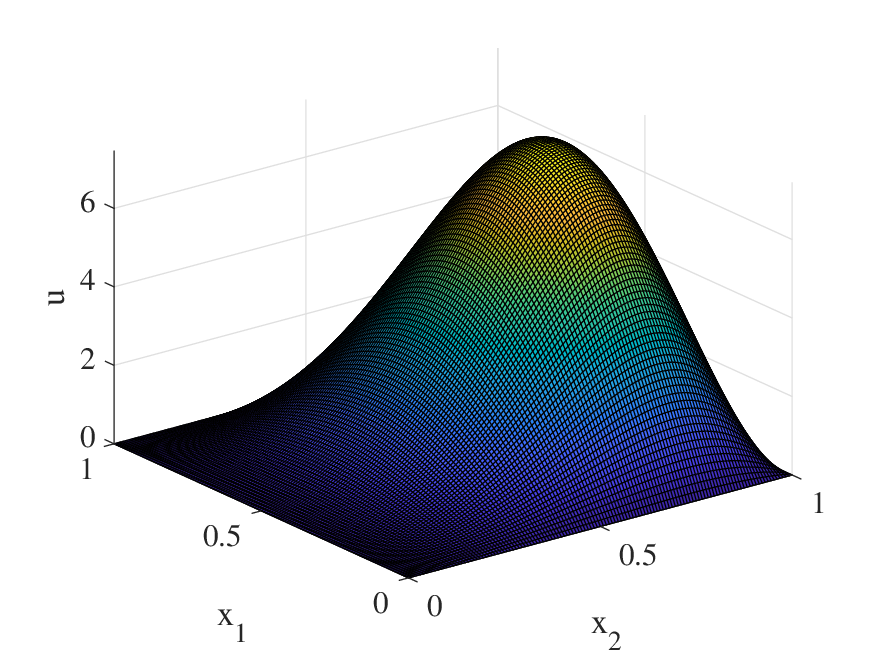}
        \caption{Numerical control $u$}
        \label{fig:numerical-control}
    \end{subfigure}
    \vspace{5pt} 
    \begin{subfigure}{0.4\linewidth}
        \centering
        \includegraphics[width=\textwidth]{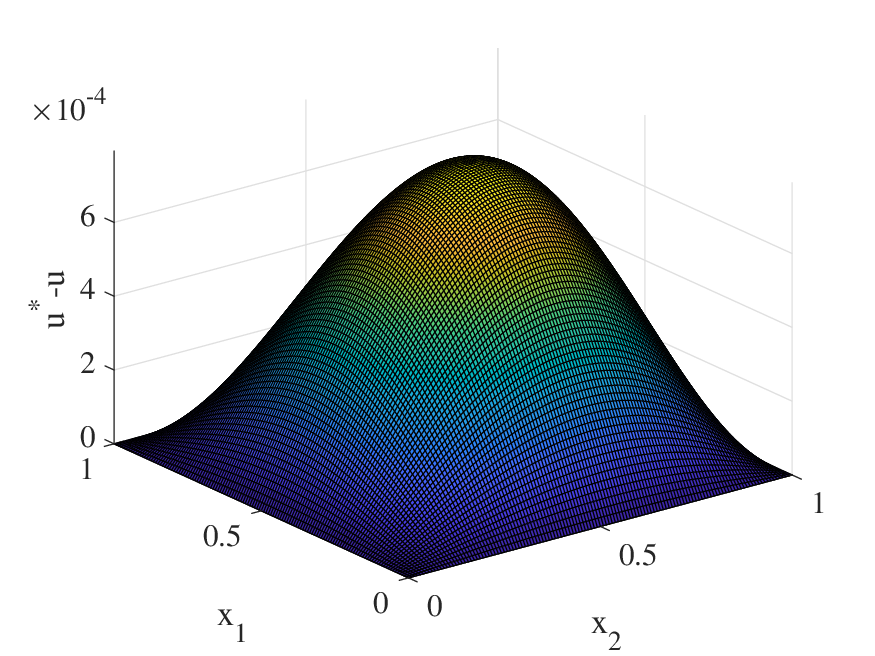} 
        \caption{Error:$u^*-u$}
        \label{fig:additional-result}
    \end{subfigure}

    \caption{Plots of ISSNG-L results}
    \label{Figure 1}
\end{figure}

In Figure \ref{Figure 1}, (a), (b), and (c) respectively display the exact control $u^*$, the numerical control $u$ obtained via the inexact semismooth Newton-GMRES method with nonmonotonic line search (ISSNG-L), and the difference between the exact and numerical control under a grid size of $h=\frac{1}{128}$. The visual similarity between (a) and (b) indicates that the two solutions are indistinguishable to the naked eye.

We modify the parameter $c_1$ in nonmonotonic line search: 
\begin{equation*}
   Q(z_k + \delta_kd_k)\leq\underset{1\leq j\leq k}{max}Q(z_j)+c_1\delta_k\nabla Q(z_k)^Td_k 
\end{equation*}of ISSNG-L by setting it to 0.5, 1, 1.3, 1.8 and 2.3, then apply the inexact semismooth Newton-GMRES method (ISSNG) to solve Example 1 with $h=\frac{1}{64}$, with the step size fixed at 1.
Table \ref{Table 1}-\ref{Table 3} present the performance of ISSNG-L with different values of the parameter $c_1$ and ISSNG under various grid step-sizes. The findings suggest that although convergence slows when $c_1>0.5$, the precision of $\Vert r_y\Vert $ and $\Vert r_p \Vert $ is enhanced.  However, it is not the case that larger values of $c_1$ always lead to better performance. As illustrated in the figure, when $c_1 =2.3$, the number of iterations increases, but the precision does not improve. Moreover, at $c_1=1$, the precision is the highest and the convergence is only one iteration slower than ISSNG, making it the optimal value that balances precision and convergence speed.

\begin{figure}[ht]
    \centering
    \begin{subfigure}{0.49\linewidth}
        \centering
        \includegraphics[width=1\textwidth]{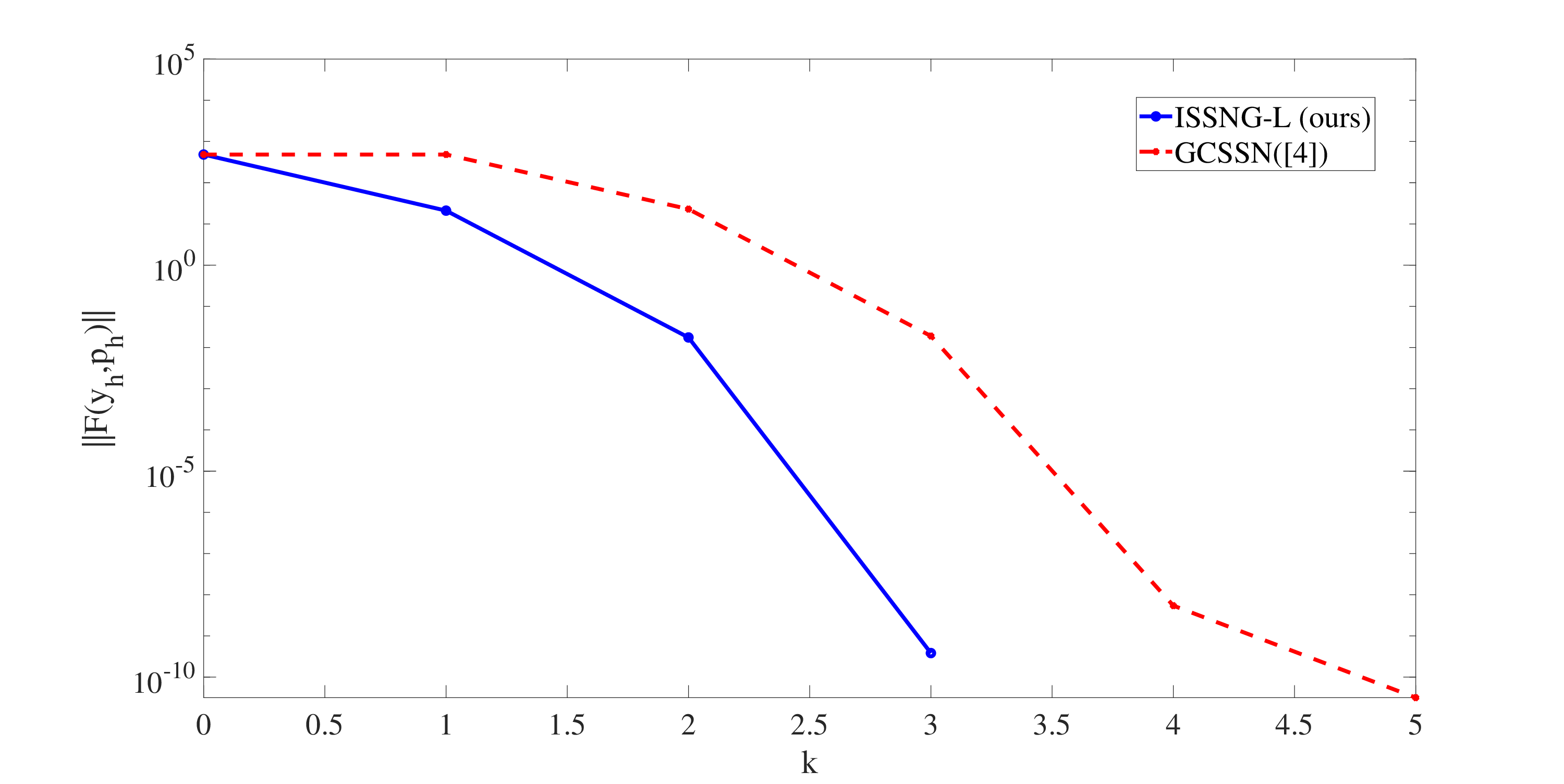} 
        \caption{$\alpha=10^{-3}$ }
    \end{subfigure}
\hspace{0.00000000000000000001\linewidth}
    \begin{subfigure}{0.49\linewidth}
        \centering
        \includegraphics[width=1\textwidth]{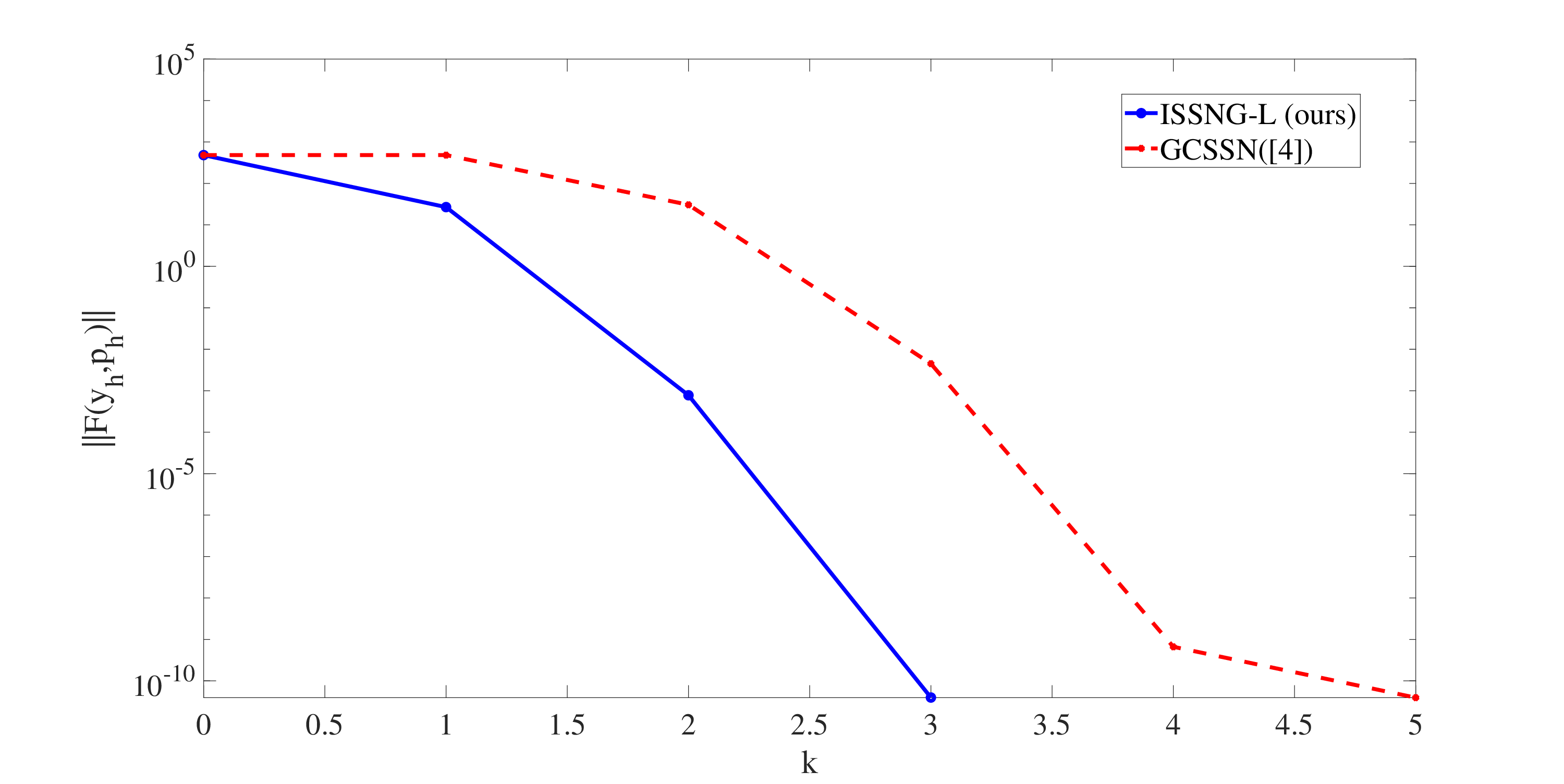} 
        \caption{$\alpha=10^{-6}$}
    \end{subfigure}
    \caption{The variations of $\Vert F(z_k)\Vert $ under the ISSNG-L with $c_1=0.5$ and GCSSN with different values of $\alpha$}
    \label{Figure 2}
\end{figure}

We conduct a comparative analysis between ISSNG-L and the GCSSN method presented in \cite{MR4447697}. Figure \ref{Figure 2}(a) and (b) respectively depict the evolution of $\Vert F(z_k)\Vert $ under the ISSNG-L and GCSSN for $\alpha=10^{-3}$ and $\alpha=10^{-6}$ with a grid step-size of $h=\frac{1}{64}$. The results indicate that ISSNG-L exhibits superior convergence performance, requiring fewer iterations to achieve the target precision. For $\alpha=10^{-6}$, Table \ref{Table 4} illustrates the values of $\Vert r_y\Vert $ and $\Vert r_p\Vert $ for both methods across various grid step-sizes. The peak memory usage, defined as the maximum memory allocated during the execution, reflects the worst-case memory requirement and is essential for judging the algorithm's scalability and practical applicability. The results indicate that when ISSNG-L and GCSSN achieve the same level of precision, ISSNG converges faster and consumes significantly less memory.


\begin{table}[H]
\centering
\caption{Numerical results of ISSNG with and ISSNG-L with $c_1=0.5$.}
\setlength{\tabcolsep}{6mm}{
\begin{tabular}{ccccccc}
\toprule
 \multicolumn{3}{c}{ISSNG} & \multicolumn{3}{c}{ISSNG-L with $c_1=0.5$} \\
\cmidrule(lr){2-4} \cmidrule(lr){5-7}
$h$& $\|r_y\|$ & $\|r_p\|$ & Iter & $\|r_y\|$ & $\|r_p\|$ & Iter \\
\midrule
$1/32$ & $1.9312e-10$ & $5.5395e-11$ & $3$ & $1.9312e-10$ & $5.5395e-11$ & $3$ \\
$1/64$ & $3.8122e-10$ & $1.0993e-10$ & $3$ & $3.8122e-10$ & $1.0993e-10$ & $3$ \\
$1/128$ & $7.4239e-10$ & $2.2824e-10$ & $3$ & $7.4239e-10$ & $2.2824e-10$ & $3$ \\
\bottomrule
\end{tabular}}
\label{Table 1}
\end{table}

\begin{table}[H]
\centering
\caption{Numerical results of ISSNG-L with $c_1=1$ and ISSNG-L with $c_1=1.3$.}
\setlength{\tabcolsep}{6mm}{
\begin{tabular}{ccccccc}
\toprule
 \multicolumn{3}{c}{ISSNG-L with $c_1=1$} & \multicolumn{3}{c}{ISSNG-L with $c_1=1.3$} \\
\cmidrule(lr){2-4} \cmidrule(lr){5-7}
$h$& $\|r_y\|$ & $\|r_p\|$ & Iter & $\|r_y\|$ & $\|r_p\|$ & Iter \\
\midrule
$1/32$ & $4.0199e-12$ & $2.8421e-14$ & $4$ &$ 4.3386e-12$ & $2.7785e-13$ & $4$ \\
$1/64$ & $3.2245e-11$ & $2.1776e-13$ & $4$ & $3.1712e-11$ & $5.8716e-13$ & $4$ \\
$1/128$ & $2.5026e-10$ & $1.7276e-12$ & $4$ & $2.5508e-10$ & $2.0150e-12$ & $4$ \\
\bottomrule
\end{tabular}}
\label{Table 2}
\end{table}

\begin{table}[H]
\centering
\caption{Numerical results of ISSNG-L with $c_1=1.8$ and ISSNG-L with $c_1=2.3$.}
\setlength{\tabcolsep}{6mm}{
\begin{tabular}{ccccccc}
\toprule
 \multicolumn{3}{c}{ISSNG-L with $c_1=1.8$} & \multicolumn{3}{c}{ISSNG-L with $c_1=2.3$} \\
\cmidrule(lr){2-4} \cmidrule(lr){5-7}
$h$& $\|r_y\|$ & $\|r_p\|$ & Iter & $\|r_y\|$ & $\|r_p\|$ & Iter \\
\midrule
$1/32$ & $4.0832e-12$ & $2.6091e-14$ & $5$ & $6.0902e-12$ & $1.5501e-12$ & $22$ \\
$1/64$ & $3.1313e-11$ & $2.1123e-13$ & $5$ & $3.2380e-11$ & $ 2.1820e-13$ & $23$ \\
$1/128$ & $2.4876e-10$ & $1.7113e-12$ & $5$ & $2.5180e-10$ & $1.7185e-12$ & $22$ \\
\bottomrule
\end{tabular}}
\label{Table 3}
\end{table}

\begin{table}[!h]
\centering
\caption{Numerical results of ISSNG-L and GCSSN.}
\setlength{\tabcolsep}{1.5mm}{
\begin{tabular}{ccccccccccc}
\toprule
 \multicolumn{5}{c}{ISSNG-L} & \multicolumn{5}{c}{GCSSN} \\
\cmidrule(lr){2-6} \cmidrule(lr){7-11}
$h$ & $\|r_y\|$ & $\|r_p\|$ & Iter & Peak memory &CPU& $\|r_y\|$ & $\|r_p\|$ & Iter & Peak memory & CPU\\
\midrule
$1/32$ & $4.8897e-12$ & $1.6873e-15$ & $3$ & $2.3142$ MB & 0.39s&$4.5280e-12$ & $1.4498e-15$ & $5$ &$9.4167$ MB& 0.87s\\
$1/64$ & $3.9036e-11$ & $3.5435e-15$ & $3$ & $9.5519$ MB &8.11s&$4.0065e-11$ & $3.5767e-15$ & $5$ &$142.6964$ MB& 19.65s\\
$1/128$ & $2.5026e-10$ & $1.7276e-12$ & $3$ & $38.7891$ MB& 56.85s&$2.5697e-10$ & $1.7131e-12$ & $5$ &$2218.8765$ MB&632.49s \\
\bottomrule
\end{tabular}}
\label{Table 4}
\end{table}

\subsection{Example 2}
Let$S(y)= y^3+y$, $u_a=-\infty$, $u_b=\infty$, $f=0$, $\alpha=10^{-3}$,
\begin{equation*}
    y_d=sin(2\pi x_1)sin(2\pi x_2)e^{2x_1}/6.
\end{equation*}
In this case, we make some certain modifications building upon Example 1. 

\begin{figure}[ht]
  \centering
  \includegraphics[width=0.8\textwidth]{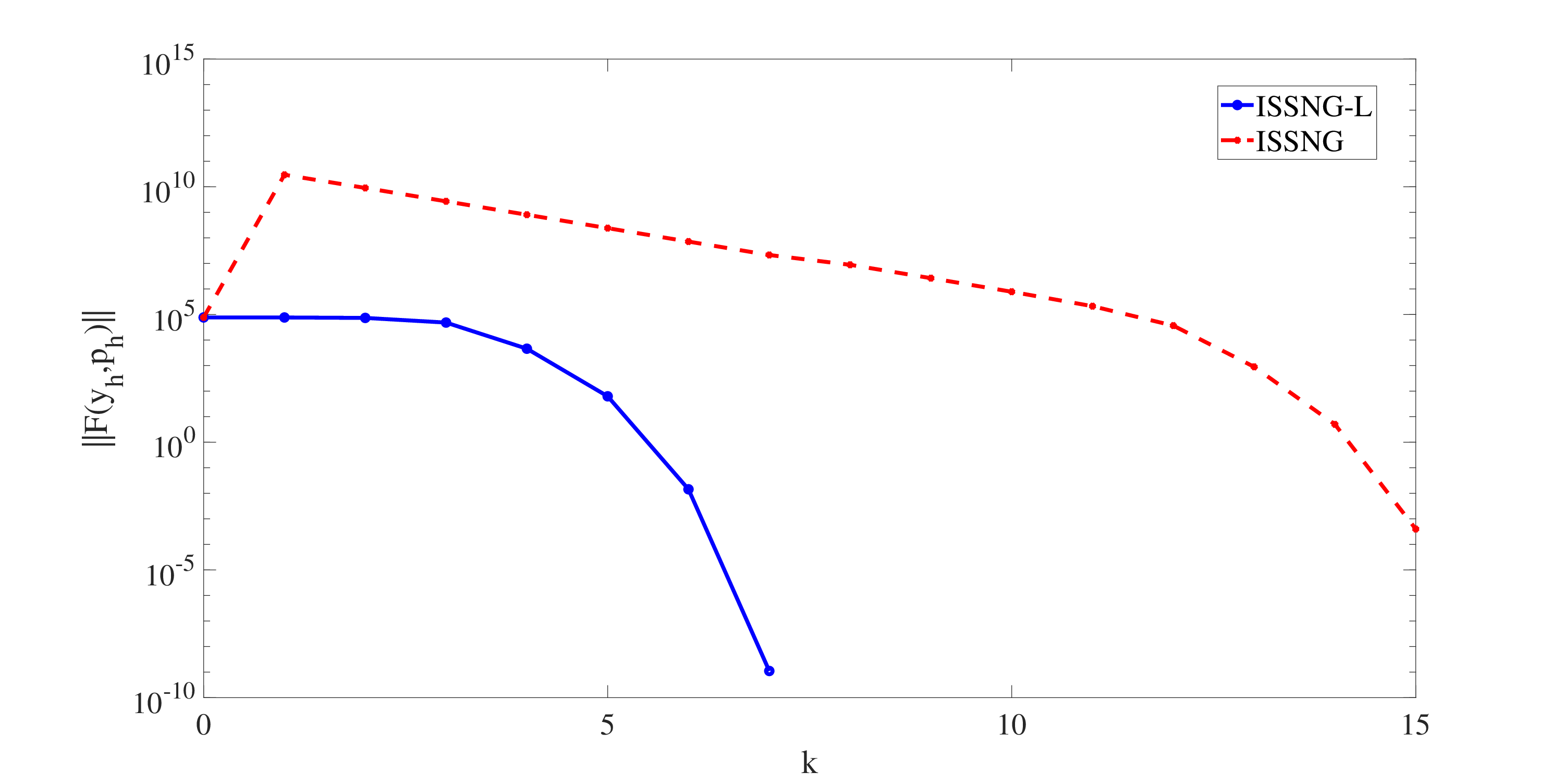} 
  \caption{The variations of $\Vert F(z_k)\Vert$ under the ISSNG-L and ISSNG}
  \label{Figure 3}
\end{figure}

The primary objective of constructing this example is to demonstrate the superiority of ISSNG-L over ISSNG for complex problems. Figure \ref{Figure 3} illustrates the evolution of $\Vert F(z_k)\Vert $ for ISSNG-L and ISSNG with a grid step-size of $h=\frac{1}{64}$ and an initial guess of matric filled with zeros. The result indicates that ISSNG-L consistently identifies a rapidly descending direction, even when the exact solution is unknown. Moreover, ISSNG-L remains unaffected by the initial guess and achieves faster convergence compared to ISSNG. Table \ref{Table 5} depicts the convergence behavior of ISSNG-L and ISSNG when the initial guesses are set to be matrices filled entirely with \textbf{0}, \textbf{1}, and \textbf{2}, respectively, across different grid step-sizes.  It should be noted that these specific values were chosen for convenience; however, the convergence of ISSNG-L is guaranteed from an arbitrary starting point by its proven global convergence property established in Section 4 of this paper. From the reported CPU times, it is evident that the ISSNG-L requires progressively less time than GCSSN to achieve convergence with increasing grid refinement. The results corroborate the aforementioned conclusions.

\begin{table}[h]
    \centering
    \caption {Numerical results of ISSNG-L and ISSNG.}
    \setlength{\tabcolsep}{2.5mm}{
    \begin{tabular}{cccc}
    \toprule
     \multirow{2}{*}{h} & \multirow{2}{*}{Initial guess} & \multicolumn{1}{c}{ISSNG-L} & \multicolumn{1}{c}{ISSNG} \\
    \cmidrule(lr){3-4}
    & & Iter &  Iter  \\
    \midrule
    \multirow{3}{*}{$\frac{1}{32}$} 
    & \textbf{0} & 7 & 13 \\
    & \textbf{1} & 5 & 8 \\
    & \textbf{2} & 5 & 6 \\
    \midrule
    \multirow{3}{*}{$\frac{1}{64}$}
    & \textbf{0} & 7 & 15 \\
    & \textbf{1} & 5 & 11 \\
    & \textbf{2} & 5 & 8 \\
    \midrule
    \multirow{3}{*}{$\frac{1}{128}$}
   & \textbf{0} & 7 & 18 \\
    & \textbf{1} & 6 & 14 \\
    & \textbf{2} & 6 & 11 \\
    \bottomrule
    \end{tabular}}
    \label{Table 5}
\end{table}

\section{Conclusion}
We propose the inexact semismooth Newton-GMRES method with nonmonotonic line search (ISSNG-L) for solving semilinear elliptic optimal control problems with box constraints. This method combines the semismooth Newton approach with Krylov subspace methods, effectively circumventing the need for direct inversion of large-scale matrices and thereby significantly reducing memory consumption during computation. By introducing a merit function to assess the quality of the solution and utilizing it to establish a nonmonotonic line search, we ensure the global convergence of the algorithm. Moreover, we prove that the algorithm achieves local superlinear convergence as the residual control parameter $\eta_k$ approaches 0. These properties are not only theoretically established but also empirically validated through numerical experiments. The method provides an efficient tool for tackling complex PDE optimization problems under distributed control constraints, particularly in scenarios requiring high-precision solutions with limited computational resources. Its low-memory characteristics make it highly promising for applications in engineering optimization and other related fields.

\bibliographystyle{elsarticle-num}
\bibliography{ISSNGL}

\end{document}